\documentclass{bmcart}

\usepackage[utf8]{inputenc} 

\usepackage{amssymb}
\usepackage{graphicx}
\usepackage{epstopdf}
\usepackage[caption=false]{subfig}
\usepackage{color}

\usepackage{diagbox}                    
\usepackage{multirow}

\usepackage{ntheorem}

\usepackage[numbers,sort&compress]{natbib}
\usepackage[misc]{ifsym}

\usepackage{lineno}

\startlocaldefs

\newtheorem*{lemma}{Lemma}[section]
\newtheorem{definition}{Definition}[section]
\newtheorem{remark}{Remark}[section]
\newtheorem{example}{Example}[section]
\newtheorem{criterion}{Criterion}[section]
\endlocaldefs

\begin{document}

\begin{frontmatter}

\begin{fmbox}
\dochead{Research}


\title{Does the smooth planar dynamical system with one arbitrary limit cycle always exists smooth Lyapunov function?}


\author[
   addressref={aff1,aff2},                   
   corref={aff1},                       
   email={ganxiaoliang@i.shu.edu.cn.com}   
]{\inits{JE}\fnm{Xiaoliang} \snm{Gan}}
\author[
   addressref={aff1,aff2},
   email={wanghaoyu042@163.com}
]{\inits{JRS}\fnm{Haoyu} \snm{Wang}}
\author[
   addressref={aff2,aff3},
   email={aoping@sjtu.edu.cn}
]{\inits{JRS}\fnm{Ping} \snm{Ao}}
\author[
   addressref={aff2,aff3},
   email={kl47895@qq.com}
]{\inits{JRS}\fnm{Yuankai} \snm{Cao}}


\address[id=aff1]{
  \orgname{Department of Mathematics, Shanghai University}, 
  \street{ShangDa Road 99},                     %
  \postcode{200444}                                
  \city{Shanghai},                              
  \cny{China}                                    
}
\address[id=aff2]{%
  \orgname{Shanghai Center for Quantitative Life Sciences, Shanghai University},
  \street{ShangDa Road 99},                     %
  \postcode{200444}                                
  \city{Shanghai},                              
  \cny{China}                                    
}
\address[id=aff3]{
  \orgname{Department of Physics, Shanghai University}, 
  \street{ShangDa Road 99},                     %
  \postcode{200444}                                
  \city{Shanghai},                              
  \cny{China}                                    
}


\end{fmbox}


\begin{abstractbox}

\begin{abstract} 
%
A rigorous proof of a theorem on the coexistence of smooth Lyapunov function and smooth planar dynamical system with one arbitrary limit cycle is given, combining with a novel decomposition of the dynamical system from the perspective of mechanics.
We base on this dynamic structure
incorporating several efforts of this dynamic structure on fixed points, limit cycles and chaos, as well as on relevant known results, such as Schoenflies theorem, Riemann mapping theorem, boundary correspondence theorem and differential geometry theory, to prove this coexistence. We divide our procedure into three steps. We first introduce a new definition of Lyapunov function for these three types of attractors. Next, we prove a lemma that arbitrary simple closed curve in plane is diffeomorphic to the unit circle. Then, the strict construction of smooth Lyapunov function of the system with circle as limit cycle is given by the definition of a potential function. And then, a theorem is hence obtained: The smooth Lyapunov function always exists for the smooth planar dynamical system with one arbitrary limit cycle. Finally, by discussing the two criteria for system dissipation(divergence and dissipation power), we find they are not equal, and explain the meaning of dissipation in an infinitely repeated motion of limit cycle.
\end{abstract}




\end{abstractbox}
%

\end{frontmatter}



\section{Introduction}
In the study of nonlinear dynamics
\begin{equation}\label{1b1}
\dot x = f(x),
\end{equation}
where $f:{\mathbb{R}^n} \to {\mathbb{R}^n}$ and $f(x)$ is smooth, the limit cycle system is one of the archetypes which have been fascinating mathematicians \cite{Poincare-1881-p375,Hilbert-1902-p437,Qin-1984-p,Ye-1986-p,TianHanXu-2019-p1561}, biologists \cite{ZehringWheelerReddyKonopka-1984-p369,Yuan-2013-p62109,Tang-2013-p12708}, physicists \cite{EmelianovaKuznetsovTurukina-2014-p153} and engineers \cite{ShafiArcakJovanovicPackard-2013-p3613} for over 100 years.
Different approaches have been developed to the study of qualitative analysis of system (\ref{1b1}), the Lyapunov function has been one of the most efficient approaches, which, however, depends on its existence.
In view of the important position of limit cycle system in dynamical system research,  proving its existence of smooth Lyapunov function is not only necessary in theory, but also has important engineering application value.

Does the smooth planar dynamical system with one arbitrary limit cycle always exists smooth Lyapunov function?
Although many researchers have paid attention to this issue, there still no definite positive result had been reached.
Such as, Hirsch, Smale and Devaney have noticed it and given a negative statement: "If $L$ is a strict Lyapunov function for a planar system, then there are no limit cycles\cite{HirschSmaleDevaney-2013-p193}".
A systematical study of qualitative analysis of dynamical system from the complete Lyapunov function of Conley theory basically began with the paper\cite{Conley-1978-p39}, about which some efforts \cite{Hurley-1998-p245,FarberKappelerLatschevZehnder-2004-p1451,Pageault-2009-p2426,Souza-2012-p322,ArgaezGieslHafstein-2018-p1,ArgaezGieslHafstein-2019-p1622} have been made.
Here, the complete Lyapunov functions is continuous for flows defined on compact metric spaces, which is a constant real function on the orbit of chain recurrent set and strictly decreases along all other orbits.
Nevertheless, the smooth Lyapunov function should exists its Lie derivative along the trajectory.
What's more, the potential function(here it is called Lyapunov function) is also taken into consideration by researchers in some practical applications \cite{Wang-2009-p3038,Ge-2012-p23140,Vincent-2018-p62203}, while they are special cases.
Then, is there a potential function method that can systematically study limit cycle systems?
The answer is yes.
Ao et al. \cite{Ao-2004-p25,Yuan-2012-p7010} in 2004 divided the dynamical system (\ref{1b1}) into three parts from the perspective of mechanics  to systematically study the behavior of the dynamic system
\begin{equation}\label{1j5000}
(S(x) + T(x))\dot x =  - \nabla \phi (x),
\end{equation}
which is  assumed to have another form
\begin{equation}\label{1j507}
\dot x =  - (D(x) + Q(x))\nabla \phi (x),
\end{equation}
here,  the single-valued scalar function $\phi (x)$ is potential function which has been proved equivalent to Lyapunov function in \cite{Yuan-2013-p10505}, the symmetric semi-positive matrix $S(x)$ corresponds to friction, the transverse matrix $T(x)$ corresponds to Lorentz force, the symmetric semi-positive matrix $D(x)$ is the diffusion matrix which can be selected by different diffusion modes,
the anti-symmetric matrix $Q(x)$ denotes Poisson bracket, and they satisfy $S(x) + T(x) = {(D(x) + Q(x))^{ - 1}}$.

Based on this novel dynamic structure, Ao et al. have made some researches on Lyapunov functions for limit cycle systems.
For instance,
in 2006, Zhu et al.  first explicitly constructed a global smooth Lyapunov function in a limit cycle system \cite{Zhu-2006-p817}.
By using the geometric method, the Lyapunov function for a piecewise linear system with limit cycle is constructed by Ma et al. in 2013 \cite{Ma-2013-p}.
Based on  A-type stochastic integral, the Lyapunov function for a class of nonlinear system, Van der Pol type system, is studied in 2013 \cite{Yuan-2013-p62109}.
In 2013, Tang et al. determined the dynamical behavior of the competitive Lotka-Volterra system by Lyapunov function \cite{Tang-2013-p12708}.

Inspired by the structural method in \cite{Zhu-2006-p817} and the strict results of special cases in \cite{Ma-2013-p,Yuan-2013-p62109,Tang-2013-p12708}, we will continue to base on the dynamic structure of Ao, and further give a definite positive result to the coexistence of smooth Lyapunov function and smooth planar dynamical system with one arbitrary limit cycle in this paper.

The rest of this paper is organized as follows.
In the next section, we combine with the efforts on fixed point, limit cycle and chaos of Ao et al., introduce a definition of Lyapunov function for these three types of attractors, and show our main results.
Some numerical examples will be presented in Section 3 to illustrate our results.
In Section 4, we notice one contradictory phenomenon about the existence of Lyapunov function for limit cycle system, and discuss and analyze it by the dynamic structure of Ao.
We conclude with Section 5.

\section{Problem formulation}
In this section, a definition of Lyapunov function for the three types of attractors, fixed point, limit cycle and chaos, is given, and the main results are shown.

\subsection{The generalized definition of Lyapunov function}
Wolfram \cite{Wolfram-2002-p961} gave a geometrical classification of attractors:
Fixed point, limit cycle and chaos.
The fixed point corresponds to zero-dimensional subset of state space, limit cycle to one-dimensional subset. Strange attractor often has a nested structure with non-integer fractal dimension.
The definition of Lyapunov function just for fixed point is clearly given, such as \cite{Wiggins-2003-p22,Hsu-2013-p140,HirschSmaleDevaney-2013-p193,Alligood-1996-p305}.
Therefore, the generalization of Lyapunov function is very necessary.

By going back to original motivation of Lyapunov by energy function, we take the efforts on these three types of attractors of Ao et al.\cite{Kwon-2005-p13029,Zhu-2006-p817,MaTanYuanYuanAo-2014-p1450015} into consideration, and give a generalized definition of Lyapunov function, which not only retains the idea that Lyapunov function does not increase along the trajectory in the usual definition, but also has a wider scope of application.

\begin{definition}\label{def1}
For a smooth autonomous system $\dot x = f(x),$ $x\in \mathbb{R}^n$, let ${x^ {**} }$ be one type of the attractors above. And $\phi: \mathbb{R}^n\rightarrow \mathbb{R}$ is a continuous differentiable function.  If it satisfies
\begin{description}
  \item[(i)] for all $x\in \mathbb{R}^n$, $\phi (x)$ has an infimum .
  \item[(ii)] for all $x\in \mathbb{R}^n$, $\phi(x)$ does not increase along the trajectory, that is
  \begin{equation}\label{1j1017}
  \dot \phi (x) = \frac{{d\phi }}{{dt}} \le 0.
  \end{equation}
\end{description}
Then, $\phi(x)$ is said to be the Lyapunov function for the three types of attractors: Fixed point, limit cycle and chaos.
\end{definition}

\begin{remark}
The above definition extends the Lyapunov function to the three types of attractors.
The condition (i) does not require the positive definiteness of Lyapunov function, and it extends the condition of general definition;
The condition (ii) ensures that it doesn't increase along trajectory.
\end{remark}

\subsection{Main results}
In this part, just the smooth planar dynamical systems are considered.

\begin{figure}[h]
\centering
\includegraphics[width=4in]{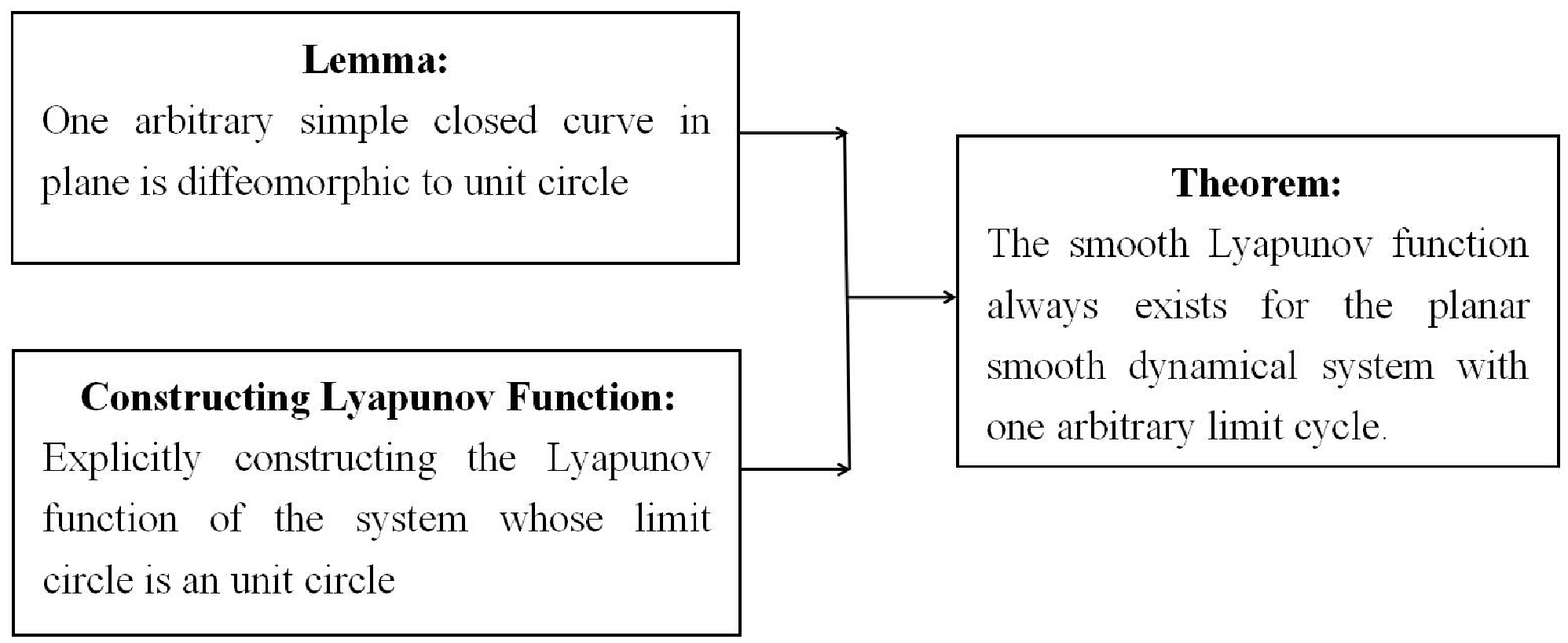}
\caption{The flow chart of main results}
\label{flowdiagram}
\end{figure}

The Fig. \ref{flowdiagram} shows our thoughts and results.
Here, main contents is briefly introduced in the following:
Firstly, inspired by Schoenflies theorem \cite{Rolfsen-2003-p9}, we combine with Riemann mapping theorem\cite{Ahlfors-1979-p135}, boundary correspondence theorem \cite{Goluzin-1969-p31} and the differential geometry theory to discuss the limit cycle which can be expressed as simple closed curves in the complex plane, and obtain a lemma that one arbitrary simple closed curve in plane is diffeomorphic to unit circle.
Secondly, we base on the thought of potential energy\cite{Strogatz-2014-p28}, and construct the smooth Lyapunov function for the system whose limit circle is an unit circle.
Thirdly, using the results above, we can summarize a new and important theorem of dynamical systems: The smooth Lyapunov functions always exist for the planar smooth dynamical systems with one arbitrary limit cycle.

\subsubsection{To prove one arbitrary simple closed curve in plane diffeomorphic to unit circle.}
\begin{lemma}
One arbitrary simple closed curve in plane is diffeomorphic to unit circle.
\end{lemma}

\emph{Proof}:
Although the result looks intuitive, it is difficult to prove, which is a tough task in complex function and topology. Since the content is too long, we divide the proof into three parts and just give the main idea here. Please refer to Appendix A for details.
\begin{description}
  \item[Part I:] Prove there exists a bijection mapping the region bounded by a simple closed curve to the unit disk;
  \item[Part II:] Prove the obtained bijection in part I is a bijection on boundary;
  \item[Part III:] Prove the obtained bijection is a diffeomorphic mapping on boundary.
\end{description}

\subsubsection{Construct smooth Lyapunov function}
The Lemma presents that an arbitrary smooth limit cycle in plane can be transformed into unit circle by appropriate scaling and translation transformation.
So constructing the smooth Lyapunov function for the smooth system whose limit circle is an unit circle becomes another basic and key step.

For the sake of simplicity, the polar coordinate system of a plane dynamical system with the circle as the limit cycle can be expressed as
\begin{equation}\label{206}
\left\{ \begin{array}{l}
\dot r = \Upsilon_0(r)\\
\dot \theta  = \psi (r,\theta)
\end{array} \right.,
\end{equation}
here, $r = \sqrt {x^2 + y^2} $, the function $\Upsilon_0(r)$ and $ \psi (r,\theta)$ are smooth.
The (\ref{181}) of the following Example \ref{example2} is a special case of this situation.
The polar coordinate system of a plane dynamical system with the circle as the limit cycle also has the other two forms.
By the detailed discussion in Appendix B blow, we find it is appropriate to use (\ref{206}).

In this place, we just consider the radial system $\dot r = \Upsilon_0(r)$ in (\ref{206}).
Set $r^{*}>0$, and $r^{*}$ is the fixed point of the one-dimensional radial system, that is $\Upsilon_0(r=r^{*}) = 0$.
Based on the thought of potential energy\cite{Strogatz-2014-p28}, the potential $\phi(r)$ is defined as
\begin{equation}\label{1.27}
\Upsilon_0(r) =- \frac{{d\phi (r)}}{{dr}}.
\end{equation}
To verify this definition, think of $r$ as a function of $t$, and calculate the time-derivative of $\Upsilon_0(r(t))$ by using the chain rule, it has
\begin{equation}\label{1.28}
\frac{{d\phi (r)}}{{dt}} = \frac{{d\phi (r)}}{{dr}}\frac{{dr}}{{dt}},
\end{equation}
Combine with the definition of potential (\ref{1.27}), for the first order system $\frac{{dr}}{{dt}} =  - \frac{{d\phi (r)}}{{dr}}$, it can derive
\begin{eqnarray}\label{1.29}
\frac{{d\phi (r)}}{{dt}} &=& \frac{{d\phi (r)}}{{dr}}\frac{{dr}}{{dt}}   \nonumber\\
   &=& - {\left( {\frac{{d\phi (r)}}{{dr}}} \right)^2} \nonumber\\
   &\le&  0,
\end{eqnarray}
so $\phi(r)$ decreases along the trajectory.

Then, by the definition of potential (\ref{1.27}), it gets
\begin{equation}\label{207}
\phi (r) =  - \int {\Upsilon_0(r)dr} .
\end{equation}
And then, by changing the polar coordinates of (\ref{207}) into the Cartesian coordinates, the potential function in Cartesian coordinates is obtained.
As mentioned in the dynamic structure (\ref{1j5000}) of Ao, the equivalence between generalized Lyapunov function and potential function is proved in literature \cite{Yuan-2013-p10505}. In following, the potential function is called Lyapunov function.

Combining the proof of Lemma and the explicit construction of smooth Lyapunov function, a theorem is obtained:

\textbf{Theorem}~\emph{The smooth Lyapunov function always exists for the smooth planar  dynamical system with one arbitrary limit cycle.}

\section{Examples}
In this section, some examples are given to illustrate the conclusions obtained in the previous section, respectively.

Firstly, the Example \ref{example1} can verify the conclusion of Lemma.
\begin{example}\label{example1}
Vibration equation\cite{Abdelkader-1998-p308}
\begin{equation}\label{1a35a}
\ddot x + (4{x^2} - 1)\dot x + x - {x^3} + {x^5} = 0.
\end{equation}
\end{example}

\emph{Solve:}
With $x = x,\dot x = y$ we have the equivalent vector equation of (\ref{1a35a})
\begin{equation}\label{1a36a}
\left\{ \begin{array}{l}
\dot x = y\\
\dot y =  - (4{x^2} - 1)y - x + {x^3} - {x^5}
\end{array} \right..
\end{equation}
By the invertible transformation
\begin{equation}\label{1a37}
\left\{ \begin{array}{l}
x = u\\
y = v + u - {u^3}
\end{array} \right.,
\end{equation}
the (\ref{1a36a}) can be expressed as
\begin{equation}\label{1a38a}
\left\{ \begin{array}{l}
\dot u = v + u - {u^3}\\
\dot v =  - {u^2}v - u
\end{array} \right.,
\end{equation}
then, by the polar transformation
\begin{equation}\label{1a39a}
\left\{ \begin{array}{l}
u = r\cos \theta \\
v = r\sin \theta
\end{array} \right.,
\end{equation}
the (\ref{1a38a}) can be expressed as
\begin{equation}\label{1a41a}
\left\{ {\begin{array}{*{20}{l}}
{\dot r = r\left( {1 - {r^2}} \right){{\cos }^2}\theta }\\
{\dot \theta  =  - 1 - \cos \theta \sin \theta }
\end{array}} \right.,
\end{equation}
here, set ${\Upsilon _0}(r){= }r(1 - {r^2}),{\Upsilon _1}(\theta ) = {cos^2}\theta $.

Because ${\Upsilon _0}(1) \equiv 0$, ${\Upsilon _1}(\theta ) \not \equiv 0$, it can obtain that the non-zero fixed point $r=1$ of ${\Upsilon _0}(r)$ in (\ref{1a41a}) corresponds to one limit cycle, and that ${\Upsilon _1}(\theta )$  does not affect the size and position of limit cycle.
So by reversible transformation(\ref{1a37}), a general limit cycle of system (\ref{1a36a}) can be translated into the one which has the shape of unit circle of system (\ref{1a38a}).
\begin{figure}[htbp]
\centering
\subfloat[The phase diagram of (\ref{1a36a}).]{%
\resizebox*{5.9cm}{!}{\includegraphics{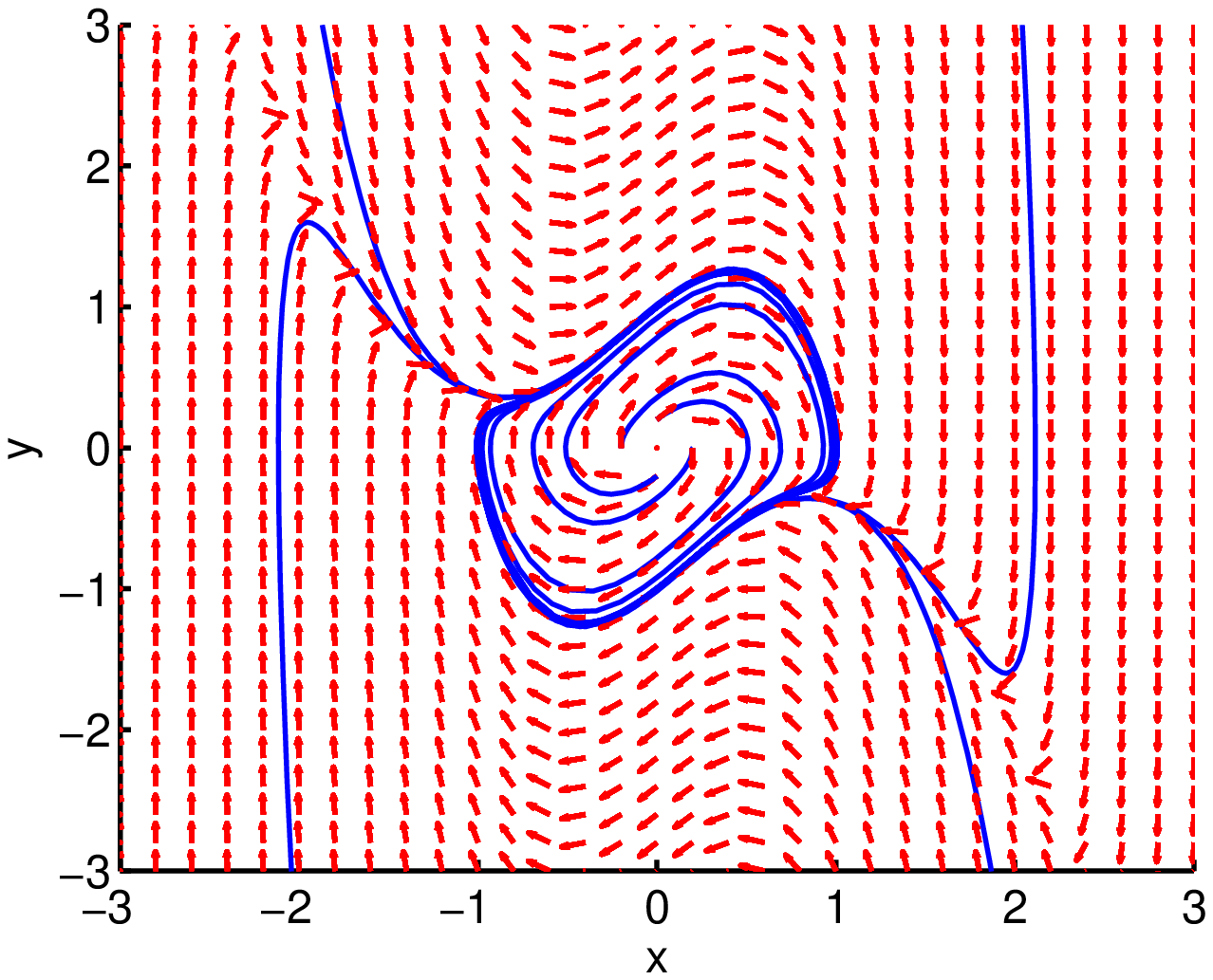}}}\hspace{5pt}
\subfloat[The phase diagram of (\ref{1a38a}).]{%
\resizebox*{5.9cm}{!}{\includegraphics{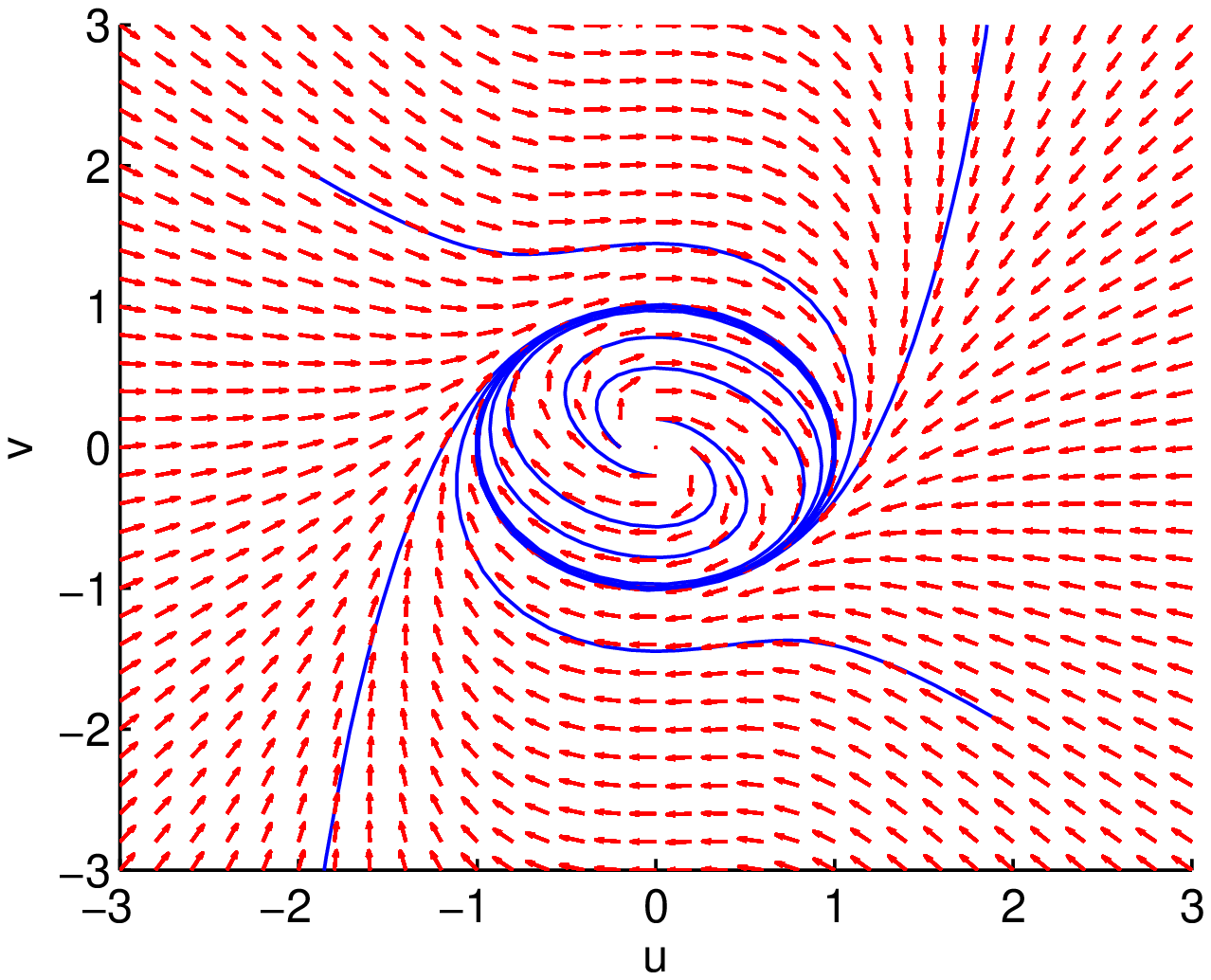}}}
\caption{The phase diagram of (\ref{1a36a}) and (\ref{1a38a}).} \label{figure4}
\end{figure}

The Fig. \ref{figure4} shows the phase diagrams of the original system (\ref{1a36a}) and transformed system (\ref{1a38a}).

Secondly, the Example \ref{example2} verifies the construction of smooth Lyapunov function with a limit cycle as unit circle.
\begin{example}\label{example2}
Find the Lyapunov function of the system\cite{Perko-2001-p195}
\begin{equation}\label{172}
\left\{ {\begin{array}{*{20}{l}}
{\dot x =  - y + x[1 - ({x^2} + {y^2})] }\\
{\dot y = x + y[1 - ({x^2} + {y^2})] }
\end{array}} \right.,
\end{equation}
\end{example}
\emph{Solve:} By polar coordinate transformation, system (\ref{172}) can be transformed into
\begin{equation}\label{181}
\left\{ \begin{array}{l}
\frac{{dr}}{{dt}} = r - {r^3},\\
\frac{{d\theta }}{{dt}} =   1,
\end{array} \right.
\end{equation}
here, $r = \sqrt {{x^2} + {y^2}} $. Obviously, system (\ref{172}) has one limit cycle $r=1$.

By (\ref{207}), the Lyapunov function of system (\ref{172}) can be derived, which has the shape of a Mongolian hat
\begin{equation}\label{182}
\phi (x,y) = \frac{1}{4}({x^2} + {y^2})({x^2} + {y^2} - 2).
\end{equation}

Then, we can verify that the Lyapunov function(\ref{182}) doesn't increase along trajectory
\begin{eqnarray}\label{183}
\frac{{d\phi (x,y)}}{{dt}}
   &=&  - ({x^2} + {y^2}){({x^2} + {y^2} - 1)^2} \le 0 .
\end{eqnarray}

\begin{figure}[htbp]
\centering
\subfloat[The phase diagram.]{%
\resizebox*{5.9cm}{!}{\includegraphics{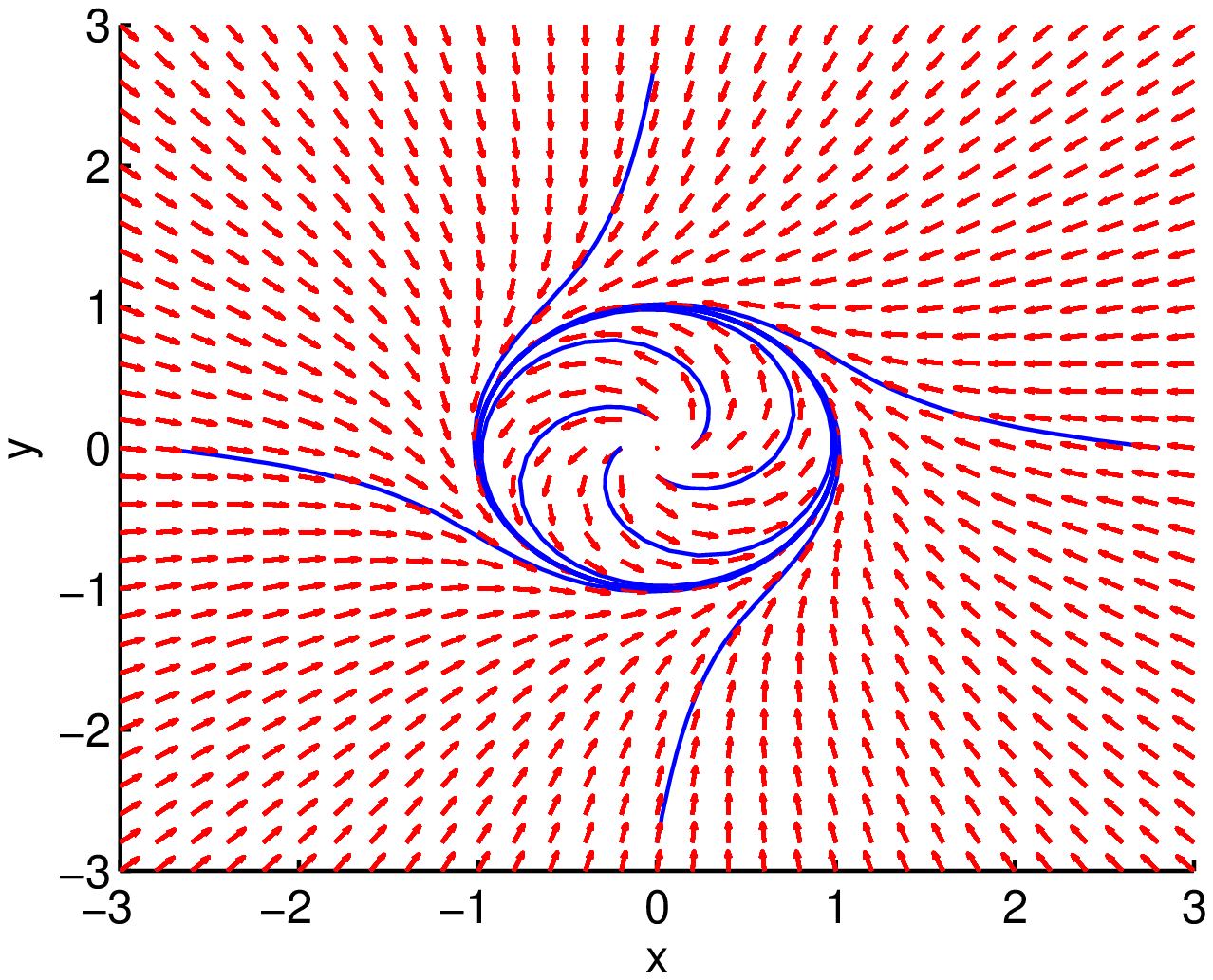}}}\hspace{5pt}
\subfloat[Lyapunov function.]{%
\resizebox*{5.9cm}{!}{\includegraphics{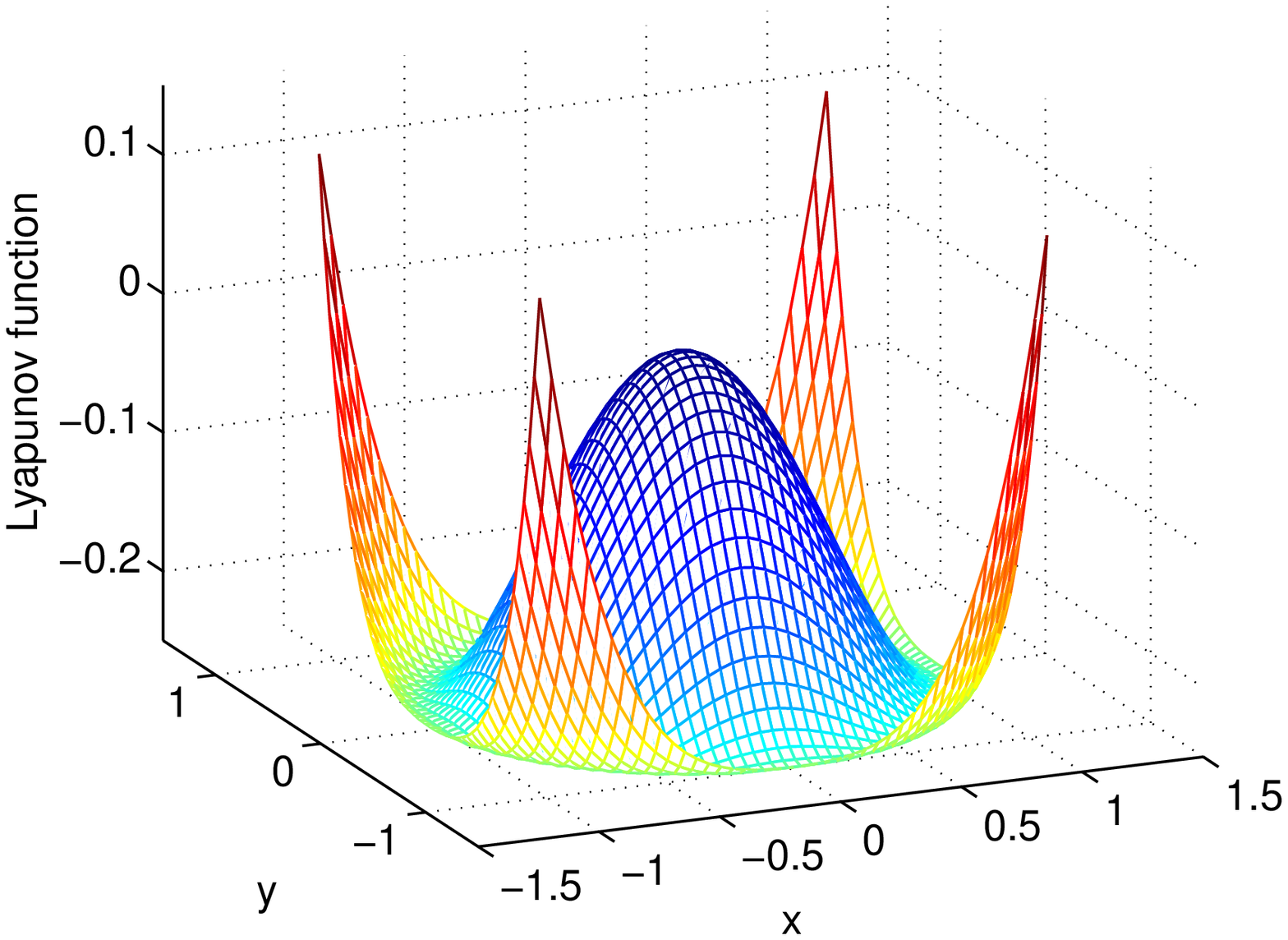}}}
\caption{The phase diagram and the Lyapunov function of system (\ref{172}).}
\label{figure5}
\end{figure}

The Fig. \ref{figure5} shows the phase diagram and Lyapunov function of system (\ref{172}).

Finally, combining the Example \ref{example1} and Example \ref{example2}, we will further derive out the Lyapunov function of Example \ref{example1}.

\begin{example}\label{example3}
Find the smooth Lyapunov function of the system (\ref{1a36a}) in Example \ref{example1}.
\end{example}
\emph{Solve:} By the inverse transformation (\ref{1a37}) and polar transformation, the original system (\ref{1a36a}) is transformed into
\begin{equation}\label{1a41ab}
\left\{ \begin{array}{l}
\frac{{dr}}{{dt}} = r(1 - {r^2}){\cos ^2}\theta \\
\frac{{d\theta }}{{dt}} =  - 1 - \cos \theta \sin \theta
\end{array} \right..
\end{equation}
Combined with (\ref{207}), it can derive the Lyapunov function of system (\ref{1a38a})
\begin{equation}\label{1a43}
\overline \phi  (r) =  - \int {r(1 - {r^2})dr}  = \frac{{{r^2}}}{4}({r^2} - 2),
\end{equation}
it is,
\begin{equation}\label{1a44}
\overline \phi (u,v) = \frac{1}{4}({u^2} + {v^2})({u^2} + {v^2} - 2).
\end{equation}
By the inverse transformation of (\ref{1a37})
\begin{equation}\label{1a45}
\left\{ \begin{array}{l}
u = x\\
v = y - x + {x^3}
\end{array} \right.,
\end{equation}
the Lyapunov function of system (\ref{1a36a}) is obtained
\begin{equation}\label{1a46}
\phi (x,y) = \frac{1}{4}[{x^2} + {(y - x + {x^3})^2}][{x^2} + {(y - x + {x^3})^2} - 2].
\end{equation}

Furthermore, we can verify that the Lyapunov function (\ref{1a46}) doesn't increase along trajectory of system (\ref{1a36a})
\begin{eqnarray}\label{1a47}
\frac{{d\phi (x,y)}}{{dt}} &=&  \frac{{d\bar \phi (r)}}{{dr}}\frac{{dr}}{{dt}}\nonumber \\
   &=& - {r^2}{(1 - {r^2})^2}{\cos ^2}\theta  \nonumber \\
   &\le& 0.
\end{eqnarray}

\begin{figure}[h]
\centering
\includegraphics[width=2.5in]{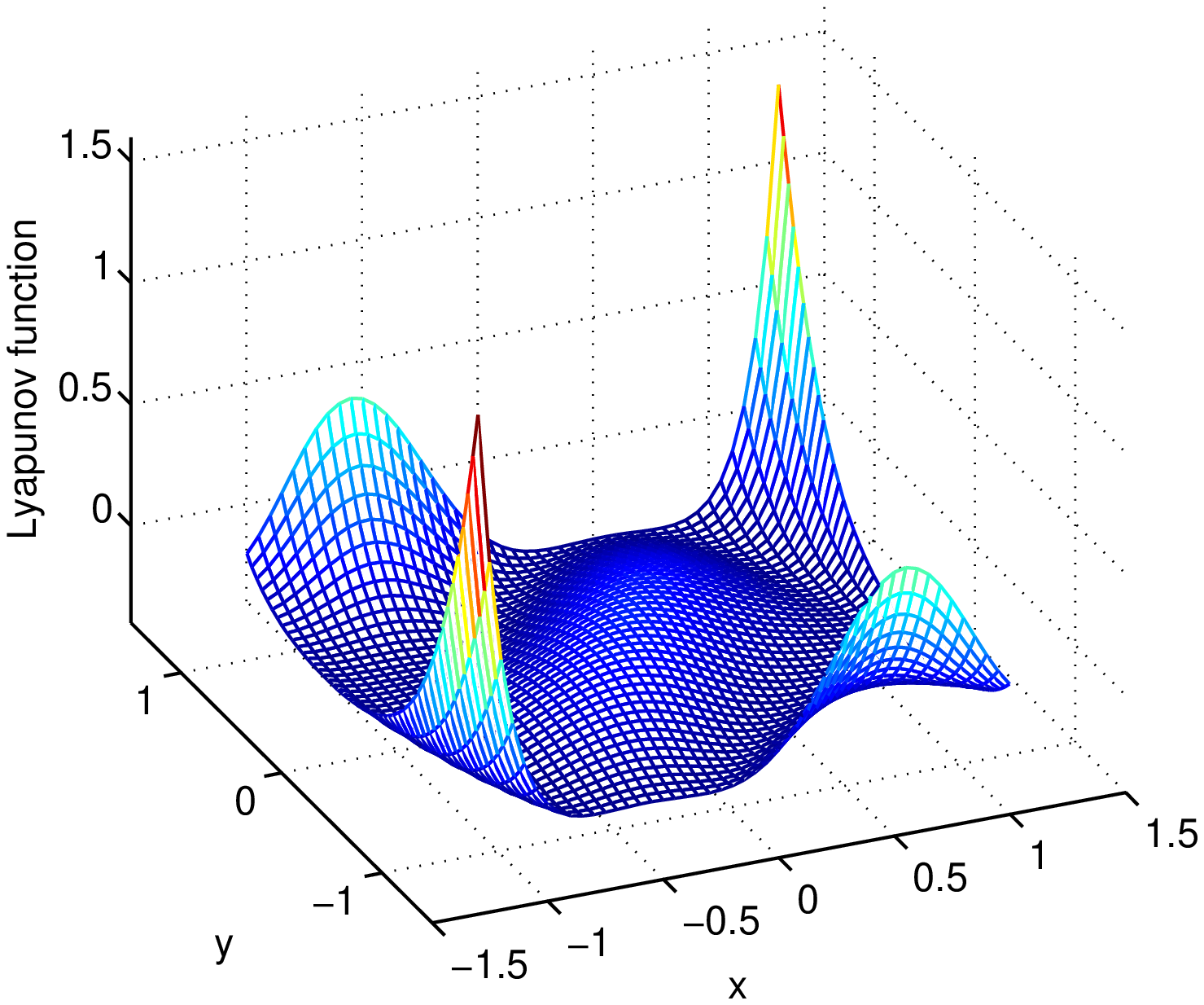}
\caption{Lyapunov function of system (\ref{1a36a})}
\label{fig13}
\end{figure}

The Figure \ref{fig13} shows the Lyapunov function of system (\ref{1a36a}).


\section{Discussion}
In system (\ref{172}) of Example \ref{example2}, the divergence on limit cycle can be obtained
\begin{eqnarray}\label{1d.67}
{\left. {divf} \right|_{{x^2} + {y^2} = 1}} &=& {\left. {\left[ {2(1 - 2{x^2} - 2{y^2})} \right]} \right|_{{x^2} + {y^2} = 1}} \nonumber \\
   &=& -2.
\end{eqnarray}
Then, a contradictory phenomenon is noticed: On the limit cycle, the system (\ref{172}) is dissipative by the divergence criterion \cite{Layek-2015-p27,Huang-2010-p35}, the trajectory moves infinitely on the limit cycle.

How to explain the meaning of dissipation in an infinitely repeated motion of a limit cycle?
To the best of our knowledge, no reasonable explanation for this phenomenon has been found yet.
Here, we will combine with the dynamic structure (\ref{1j5000}) of Ao and analyze it from the perspective of system dissipation.

There are two criteria of dissipation:

\begin{criterion}
Dissipation can be defined by dissipative power via friction force ${F_{friction}}$ \cite{Goldstein-2001-p24}:
\begin{eqnarray}\label{1j1004}
{H_P}&=& {F_{friction}} \cdot ( - \dot x) \nonumber\\
     &=& ( - S(x)\dot x) \cdot ( - \dot x)  \nonumber\\
     &=& {\dot x^\tau }S(x)\dot x   \nonumber\\
     &\geq&0,
\end{eqnarray}
here, ${F_{friction}} =  - S(x)\dot x$, the semi-positive definite symmetric matrix $S(x)$ is the friction matrix which guarantees the nonnegativity of $H_P(x)$, $x = {({x_1},{x_2}, \cdots ,{x_n})^\tau }$,  $\tau$ is transpose symbol.
This criterion classifies the dynamics into dissipative or conservative according to the change of "energy function" or "Hamiltonian".
When ${H_P}(x)>0$, the system is dissipative. When ${H_P}(x) =0$, the system is conservative. Since the $S(x)$ is a semi-positive definite symmetric matrix, $H_P<0$ is impossible.
This method is used in physics.
\end{criterion}

\begin{criterion}
Dissipation is defined as divergence \cite{Layek-2015-p27,Huang-2010-p35}:
\begin{eqnarray}\label{1j1005}
div f(x)&=& \sum {\frac{{\partial {f_i}}}{{\partial {x_i}}}}\nonumber\\
     &<&0,
\end{eqnarray}
This criterion classifies the dynamics into dissipative or conservative according to the the change of phase space volume.
When $div f(x) < 0$, the phase space volume decreases and the system is dissipative. When $div f(x) =0$, the phase space volume remains unchanged and the system is conservative.
When $div f(x) > 0$, the phase space volume increases.
This method is used in mathematics.
\end{criterion}

\begin{remark}
These two dissipation criteria are generally considered to be equivalent.  However, we find that they are not equivalent on the limit cycle.
\end{remark}

Combined with the Lyapunov function (\ref{182}) and paper \cite{Yuan-2013-p10505}, the $S(x,y)$ of system (\ref{172}) is obtained
\begin{equation}\label{1b2}
S(x,y) = \frac{{1 - {x^2} - {y^2}}}{{1 + {{(1 - {x^2} - {y^2})}^2}}}\left[ {\begin{array}{*{20}{c}}
1&0\\
0&1
\end{array}} \right],
\end{equation}
then, the corresponding dissipative power ${H_P}(x,y)$  is obtained on its limit cycle ${{x^2} + {y^2}{ = 1}}$
  \begin{eqnarray}\label{1j1007}
{\left. {{H_P}(x,y)} \right|_{{x^2} + {y^2}{= 1}}} &=&{\left. { (x^2 + y^2){{(x^2 + y^2 - 1)}^2}} \right|_{{x^2} + {y^2} = 1}}\nonumber \\
  &\equiv& 0.
\end{eqnarray}
By comparing (\ref{1d.67}) with (\ref{1j1007}), we obtain these two dissipation criteria are not equivalent on limit cycle.

What's more, the relationship between the dissipative power criterion and Lyapunov function is derived by (\ref{183}) with (\ref{1j1007})
\begin{equation}\label{1b3}
{\frac{{d\phi }}{{dt}}} = - H_p.
\end{equation}
The physical meaning of formula (\ref{1b3}) is obvious: The reduction of energy indicates dissipation, and $H_p=0$ implies there is no dissipation.
Then, we can give an explain to the meaning of dissipation in an infinitely repeated motion of a limit cycle: On the limit cycle of system (\ref{172}), the dissipative power $H_p\equiv0$  indicates that the system is conservative. So the trajectory can move infinitely on the limit cycle with no dissipation.

\section{Conclusions}
This paper has given a detailed and positive answer to the coexistence of smooth Lyapunov function and smooth planar dynamical system with one arbitrary limit cycle.

\section*{Appendix A}
\appendix
\renewcommand{\theequation}{\thesection A.\arabic{equation}}
In this appendix, we prove the Lemma. Here, we just consider the phase curve on limit cycle. Since for some $t_0 >0$, if the time $t> t_0$, the phase curve of limit cycle system will be on limit cycle.

Although the lemma seems simple, proving it is a difficult task, which should first consider the phase curve in plane into the simple curve in  complex plane, then properly combine with Schoenflies theorem \cite{Rolfsen-2003-p9}, Riemann mapping theorem\cite{Ahlfors-1979-p135}, boundary correspondence theorem \cite{Goluzin-1969-p31} and the differential geometry theory to prove it. The following is the detailed process of proof.

\emph{Proof}:
Combined with literatures\cite{Ahlfors-1979-p135,Goluzin-1969-p31}, a strict proof process is given below, which is divided into three parts:
Part I will prove that there is a bijection mapping the region bounded by a simple closed curve to the unit disk;
Part II will prove that the obtained bijection is a bijection on boundary;
Part III will prove that the obtained bijection is a diffeomorphic mapping on boundary.
\begin{description}
  \item[Part I:]
  Prove that there exists a mapping that maps the region $\Omega$ enclosed by simple closed curve $J$ to the unit disk $W$.
  The idea is to prove the existence of such a mapping by using some property of normal family\cite{Ahlfors-1979-p135}: For arbitrary point $a  \in \Omega $, there exists an unique biholomorphic mapping $f$ that maps $\Omega$ to the unit disk $W$, such that  $f(a)=0$, $f'(a)>0$. It consists of the following six steps:
\begin{description}
  \item[Step 1:] Construct the function family $\mathcal{F}$ on the bounded region $\Omega$, and make it has the following properties
  \begin{equation}\label{1d.47}
  \begin{array}{l}
{\cal F} = \left\{ {f:f~is~an~univalent~function~on~\Omega ,\left| {f(z)} \right| < 1,} \right.\\
{\kern 1pt} {\kern 1pt} {\kern 1pt} {\kern 1pt} {\kern 1pt} {\kern 1pt} {\kern 1pt} {\kern 1pt} {\kern 1pt} {\kern 1pt} {\kern 1pt} {\kern 1pt} {\kern 1pt} {\kern 1pt} {\kern 1pt} {\kern 1pt} {\kern 1pt} {\kern 1pt} {\kern 1pt} {\kern 1pt} {\kern 1pt} {\kern 1pt} {\kern 1pt} {\kern 1pt} {\kern 1pt} \left. {f(a) = 0,f'(a) > 0} \right\}
\end{array}
  \end{equation}
  Obviously, the $\mathcal{F}$ is not empty, for example: $\frac{{z - a}}{d}$ is belong to $\mathcal{F}$, $d$ is the diameter of $\Omega$.
  \item[Step 2:] Make the derivative $f'(z)$ of the selected function $f$ in $\mathcal{F}$ have the highest possible value at $a$.

  Set $\{ \Omega _n^ * \} $ as a compact subset sequence of $\Omega$, so that $\mathop  \cup \limits_{n = 1}^\infty  \Omega _n^ *  = \Omega $ and $a \in \Omega _n^ * $. And then, $\{ f(\Omega _n^ * )\} $ is a compact subset sequence of $W = \left\{ {z:\left| z \right| < 1} \right\}$. As $n\rightarrow \infty$, $\mathop  \cup \limits_{n = 1}^\infty  f(\Omega _n^ * )$ gets bigger and bigger and tries to fill the unit circle $W$. We will choose $f$ from $\mathcal{F}$, and the derivative $f'(a)$ has the highest possible value, so that we select the function "fastest spreading" at $z=a$, and $f$ has the most opportunity to satisfy $\bigcup\limits_{n = 1}^{ + \infty } {f(\Omega _n^ * )}  = W$.

   Note
  \begin{equation}\label{1d.48}
  v = \sup \left\{ {\left| {f'(a)} \right|:f \in \mathcal{F}} \right\}.
  \end{equation}
  Because $\Omega$ is an open region, there certainly exists a neighborhood $U(a,\kappa) \subset \Omega$, $\kappa>0$.
  Let $ \forall w \in W$, $\varphi (w) = \kappa w + a$, then $\psi$ maps the unit disk $W$ to $U(a,\kappa)$, and $\varphi (0) = a$.
  For $\forall f \in \mathcal{F}$, the composite function $f \circ \varphi $ maps the unit disk to the unit disk with origin to origin.
  Use the Schwarz's lemma\cite{Ahlfors-1979-p135}, it has $\left| {f'(a) \varphi '(0)} \right| \le 1$. What's more, $\varphi '(0) = r$, it gets $\left| {f'(a)} \right| \le \frac{1}{r}$ and $v <  + \infty $.
  \item[Step 3:] The mapping function is expressed as a solution of an extremum problem through some properties of the analytic function family.

  By (\ref{1d.48}), there exists a function sequence ${f_n}(z) \in \mathcal{F}(n = 1,2, \cdots )$ such that
  \begin{equation}\label{1j1064}
  \left| {{f_n}'(a)} \right| \ge v - \frac{1}{n}.
  \end{equation}
  Because $\mathcal{F}$ is uniformly bounded on $\Omega$, by the condensation principle\cite{Goluzin-1969-p31},  there exists an inner closed uniformly convergent subsequence $\left\{ {{f_{{n_j}}}(z)} \right\}$ in $\left\{ {{f_n}(z)} \right\}$. Set $\left\{ {{f_{{n_j}}}(z)} \right\}$ uniformly converges to some holomorphic function $f(z)$ on $\Omega$.
  By the Weierstrass theorem\cite{Ahlfors-1979-p135}, $\left\{ {f_{{n_j}}^{'}(z)} \right\}$ uniformly converges to $f'(z)$ and $f(z)$ is the analytic function on $\Omega$. What's more, it has $f'(a) = \mathop {\lim }\limits_{j \to  + \infty } f_{{n_j}}^{'}(a) = v > 0$, which shows that $f(z)$ is not a constant. Obviously, $\left| {f(z)} \right| < 1$. So $f(z)$ is the solution to extremum problem.
  \item[Step 4:] Prove the obtained $f(z)$ is an injection.

  That's, for $\forall {z_0} \in \Omega $, the function value $w_0=f(z_0)$ cannot be taken at the point different $z_0$. So $f(z)-w_0$ has no zero point in $\Omega \backslash \left\{ {{z_0}} \right\}$. Let ${w_j} = {f_{{n_j}}}(z_0)$, by ${f_{{n_j}}} \in \mathcal{F}$, ${f_{{n_j}}}$ is an one-to-one function, so ${f_{{n_j}}}(z) - {w_j}$ has no zero point in $\Omega \backslash \left\{ {{z_0}} \right\}$.
  Clearly, sequence $\left\{ {{f_{{n_j}}}(z) - {w_j}} \right\}$ is internal closed uniform convergence to $f(z) - {w_0}$.
  By the Hurwitz theorem\cite{Conway-1995-p152}, $f(z) - {w_0}$ has no zero point in $\Omega \backslash \left\{ {{z_0}} \right\}$, so $f(z)$ is an one-to-one function in $W$.
  \item[Step 5:] Prove the obtained $f(z)$ is a surjection.

  In other word, it is $f(\Omega)=W$.
  If $f(\Omega)\neq W$, there at least exists one point $c \in W$ so as to $c \notin f(\Omega)$, $0 < \left| c \right| < 1$.
  Combined with $f(a)=0$ and $f(z)$ being one-to-one, do the fraction linear transformation
      \begin{equation}\label{1d.49}
      {\psi _1}(z) = \frac{{f(z) - c}}{{1 - \overline c f(z)}}.
      \end{equation}
  Because ${\psi _1}(z)$ is simple connected and has no zero point, so $\sqrt {{\psi _1}(z)} $ can has the single branch in $\Omega$, which is written as ${{\psi _2}(z)}$. Obviously, it has $\left| {{\psi _2}(z)} \right| < 1$.
  Notice ${\psi _2}(a) = \sqrt { - c}  \ne 0$, do the fraction linear transformation
      \begin{equation}\label{1d.50}
      {\psi _3}(z) = {e^{i\theta }}\frac{{{\psi _2}(z) - {\psi _2}(a)}}{{1 - {{\overline \psi  }_2}(a){\psi _2}(z)}},\theta  = \arg {\psi _2}(a).
      \end{equation}
  It is easy to get $\left| {{\psi _3}(z)} \right| < 1$, $\left| {{\psi _3}(a)} \right| =0$. Because $ - c = {\psi _1}(a) = \psi _2^2(a)$, then
  \begin{eqnarray}\label{1j1016}
   \psi _3^{'}(a) &=& {e^{i\theta }}\frac{{\psi _2^{'}(a)}}{{1 - {{\left| {{\psi _2}(a)} \right|}^2}}} \nonumber \\
     &=& {e^{i\theta }}\frac{{\psi _1^{'}(a)}}{{2{\psi _2}(a)}}\frac{1}{{1 - {{\left| {{\psi _2}(a)} \right|}^2}}}  \nonumber \\
     &=& \frac{{1 + {{\left| {{\psi _2}(a)} \right|}^2}}}{{2\left| {{\psi _2}(a)} \right|}}f'(a)  \nonumber\\
     &>&  f'(a) \nonumber\\
     &=&  v.
  \end{eqnarray}
  By the property of fraction linear transformation, ${\psi _3}(z) \in \mathcal{F}$, then $\psi _3^{'}(a) \leq v$, it gets a contradiction. So $f$ is surjective.
  Combined with the step 2, it gets that $f $ is a biholomorphic mapping from $\Omega$ to $W$ and $f (a) = 0, f '(a) > 0 $.
  \item[Step 6:]Prove the $f$ is unique.

  Suppose $g$ also is a biholomorphic mapping from $\Omega$ to $W$, and $g(a) = 0,g'(a) > 0$.
  Then, $g \circ {f^{ - 1}}$ is  a mapping from the unit disk to unit disk, and maps the zero to the zero.
  By the Schwarz lemma\cite{Ahlfors-1979-p135}, it has $\left| {g({f^{ - 1}}(w))} \right| \le \left| w \right|,\forall w:\left| w \right| < 1.$ Let $w = f(z)$, then $\left| {g(z)} \right| \le \left| {f(z)} \right|$, $\forall z \in \Omega$.
  Similarly, it can get $\left| {g(z)} \right| \ge \left| {f(z)} \right|$, $\forall z \in \Omega$. Then, it has $\left| {g(z)} \right| = \left| {f(z)} \right|$, and there exists a constant $\theta_0$ satisfying $g(z) = {e^{i\theta_0 }}f(z)$.
  Because $f'(a)>0, g'(a)>0$, it gets ${e^{i\theta_0 }}=1$, and then $g(z) \equiv f(z)$.
\end{description}
  \item[Part II:] To prove that the obtained $f$ is a continuous bijection which can maps the simple closed curve $J$ to the unit circle. It consists of the following four steps:
  \begin{description}
    \item[Step 7:] Generalize the domain of $f(z)$ from $\Omega$ to its boundary $\partial \Omega$ (or written as $J$).

  Since the boundary $\partial \Omega$ is a Jordan closed curve, the points on boundary $\partial \Omega$ are accessible boundary point.
  Let $\zeta  \in \partial \Omega,\left\{ {{z_n}} \right\} \subset \Omega,\mathop {\lim }\limits_{n \to  + \infty } {z_n} = \zeta.$
  The function $f(z)$ is proved to be uniformly continuous in $\Omega$ by Yu \cite{Yu-2015-p202}.
  That has, for $\forall \varepsilon  > 0,\exists \delta  > 0,$, if $z',z'' \in \Omega$ and $\left| {z' - z''} \right| < \delta $, then
    \begin{equation}\label{1d.54}
    \left| {f(z') - f(z'')} \right| < \varepsilon .
    \end{equation}
  By Cauchy criterion for convergence,  when $n$ is large enough, there exists a positive integer $p = 1,2, \cdots $, such that $\left| {{z_{n + p}} - {z_n}} \right| < \delta $. Combine with (\ref{1d.54}), it has
     \begin{equation}\label{1d.55}
     \left| {f({z_{n + p}}) - f({z_n})} \right| < \varepsilon .
     \end{equation}
  Then, the sequence $\left\{ {f(z_n)} \right\}$ converges to a finite complex number $w$.

  Set another sequence $\left\{ {z_n^{'}} \right\} \subset \Omega$ converges to $\zeta$. Similarly, sequence $\left\{ {f(z_n^{'})} \right\}$ converges to a finite complex number $w'$. We have
     \begin{eqnarray}\label{1d.56}
      \left| {w - w'} \right| &=& \left| {w - f({z_n}) + f({z_n}) - f(z_n^{'}) + f(z_n^{'}) - w'} \right| \nonumber \\
        & \leq & \left| {w - f({z_n})} \right| + \left| {f({z_n}) - f(z_n^{'})} \right| + \left| {f(z_n^{'}) - w'} \right|.
     \end{eqnarray}
  By
     \begin{eqnarray}\label{1d.57}
     \left| {{z_n} - z_n^{'}} \right| &=& \left| {{z_n} - \zeta  + \zeta  - z_n^{'}} \right|  \nonumber\\
        &\leq&  \left| {{z_n} - \zeta } \right| + \left| {\zeta  - z_n^{'}} \right|,
     \end{eqnarray}
  and sequences $\left\{ {{z_n}} \right\},\left\{ {z_n^{'}} \right\}$ are all has the limit $\zeta$. So when the $n$ is large enough, it has
     \begin{equation}\label{1d.58}
     \left| {{z_n} - z_n^{'}} \right| < \delta .
     \end{equation}
  Combine with (\ref{1d.58}),(\ref{1d.54})and (\ref{1d.56}), it gets
     \begin{eqnarray}\label{1d.59}
     \left| {w - w'} \right| &<& \varepsilon  + \varepsilon  + \varepsilon =3\varepsilon.
     \end{eqnarray}
  For the arbitrariness of $\varepsilon$, there must be $\left| {w - w'} \right| = 0$, it is
     \begin{equation}\label{1d.57}
     {w = w'}.
     \end{equation}
  It shows
     \begin{equation}\label{1d.65}
     \exists \mathop {\lim }\limits_{z \to \zeta ,z \in G} f(z) = w, written ~as ~f(\zeta ).
     \end{equation}
    \item[Step 8:] Prove  the function $w = f(z)$ is continuous on $\partial \Omega$.

  Let $\zeta ,{\zeta ^ * } \in \partial \Omega$, $\left| {\zeta  - {\zeta ^ * }} \right| < \frac{\delta }{3}$ and $\left\{ {{z_n}} \right\},\left\{ {z_n^{\ast}} \right\}  \subset \Omega$, ${z_n} \to \zeta ,z_n^ *  \to {\zeta ^ * },n \to  + \infty $, it has
    \begin{eqnarray}\label{1d.59}
     \left| {f(\zeta ) - f({\zeta ^ * })} \right| &=& \left| {f(\zeta ) - f({z_n}) + f({z_n}) - f(z_n^ * ) + f(z_n^ * ) - f({\zeta ^ * })} \right|  \nonumber\\
       & \leq & \left| {f(\zeta ) - f({z_n})} \right| + \left| {f({z_n}) - f(z_n^ * )} \right| + \nonumber\\
       & & \left| {f(z_n^ * ) - f({\zeta ^ * })} \right| .
    \end{eqnarray}
  Due to
  \begin{eqnarray}\label{1d.60}
   \left| {{z_n} - z_n^*} \right| &=& \left| {{z_n} - \zeta  + \zeta  - {\zeta ^*} + {\zeta ^*} - z_n^*} \right| \nonumber\\
     &\leq& \left| {{z_n} - \zeta } \right| + \left| {\zeta  - {\zeta ^*}} \right| + \left| {{\zeta ^*} - z_n^*} \right|,
  \end{eqnarray}
  select the appropriate ${{z_n}, z_n ^ *} $, such that
  \begin{equation}\label{1d.61}
  \left| {{z_n} - z_n^ * } \right| < \delta .
  \end{equation}
  Then, combine with (\ref{1d.61}),(\ref{1d.54}) and (\ref{1d.59}), when $\left| {\zeta  - {\zeta ^ * }} \right| < \frac{\delta }{3}$, it has
  \begin{eqnarray}\label{1d.62}
   \left| {f(\zeta ) - f({\zeta ^ * })} \right| &<& \varepsilon  + \varepsilon  + \varepsilon  =3\varepsilon.
  \end{eqnarray}
  That is, the function $w = f(\zeta)$ is continuous on $\partial \Omega$.
    \item[Step 9:]  Prove $w = f(z)$ is an injection from $\partial \Omega$ to $\partial W = \left\{ {\left. w \right|\left| w \right| = 1} \right\}$.

        When $\zeta \in \partial \Omega$, the value of $\left| {f(\zeta )} \right| $ has two cases: $\left| {f(\zeta )} \right| < 1$ or $\left| {f(\zeta )} \right| = 1$.

  Let$\left\{ {{z_n}} \right\}\subset \Omega$, $\mathop {\lim }\limits_{n \to  + \infty } {z_n} = \zeta $, suppose $w = f(z)$ satisfies $\left| w \right| = \left| {f(z)} \right| < 1$ and there exists $z' \in \Omega$ such that $w = f(z')$.

  Select a sufficiently small neighbourhood of $z'$, written as ${\Omega_{z'}}$, when $n$ is large enough, it has ${z_n} \notin {\Omega_{z'}}$.

  According to the consistency of $f(z)$, it has $f({z_n}) \notin f({\Omega_{z'}})$. That is, $\mathop {\lim }\limits_{n \to  + \infty } f({z_n}) = f(\zeta ) = w \notin f({\Omega_{z'}})$, it is contrary to the assumption $w = f(z') \in f({\Omega_{z'}})$, so $\left| {f(\zeta )} \right| = 1$.

  Now we prove $w = f(z)$ is an injection on $\left| w \right| = 1$. That is, for  $ \forall \zeta ,{\zeta ' } \in \partial \Omega$ and $ \zeta  \neq \zeta '$, it has $f(\zeta ) \ne f(\zeta ')$.

  Let $\left\{ {{z_n}} \right\},\left\{ {z_n^{'}} \right\} \subset \Omega,\mathop {\lim }\limits_{n \to  + \infty } {z_n} = \zeta ,\mathop {\lim }\limits_{n \to  + \infty } z_n^{'} = \zeta '$, it has $\left| {f(\zeta )} \right| = \left| {f(\zeta ')} \right| = 1$.
  Select a sufficiently small neighbourhood of $\zeta$, written as ${\Omega_{\zeta}}$, when $n$ is large enough, ${z_n} \in {\Omega_\zeta },z_n^{'} \notin {\Omega_\zeta }$, then $f({z_n}) \in f({\Omega_\zeta }),f(z_n^{'}) \notin f({\Omega_\zeta })$, so $f(z_n^{'}) \to f({z_n}) \in f({\Omega_\zeta })$ is not possible. Hence, $w = f(z)$ is an injection on $\partial \Omega$.
    \item[Step 10:] Prove $w = f(z)$ is a surjection from $\partial \Omega$ to $\partial W = \left\{ {\left. w \right|\left| w \right| = 1} \right\}$.

  Without loss of generality, assume there exists ${w_0} \in \left| w \right| = 1$ such that $f(z) - {w_0} \ne 0, z \in \partial \Omega$, and it has $\left| {f(z) - {w_0}} \right| > 0$.
  Then, there must be some $r>0$, so that
   \begin{equation}\label{1d.63}
   \left| {f(z) - {w_0}} \right| > r > 0.
   \end{equation}
  It is  contradicts with the conclusion of step 8 that  $w = f(\zeta)$ is continuous on $\partial \Omega$. So $w = f(z)$ is a surjection on $\left| w \right| = 1$.
  \end{description}
  \item[Part III:] To prove that the obtained bijection $f$ is a diffeomorphic mapping. It contains the following two steps:
  \begin{description}
    \item[Step 11:] To prove that $f$ is a homeomorphism: (i) $f$ is a continuous  bijection which has been demonstrated above; (ii) $f$ and $f^{-1}$ are continuous: The continuity of $f$ has been obtained. Since the simple closed curve $J$ and the unit circle $\partial W$ are both closed sets, that is, compact, it can be known that the inverse mapping $f^{-1}$ is also continuous.

    Thus, it can obtain that $f$ is an homeomorphic mapping.
    \item[Step 12:] To prove that $f$ is differentiable.

    According to the definition of differentiable mapping \cite{Arnold-1992-p294}, Figure \ref{figDM} can be obtained.
\begin{figure}[h!]
\centering
\includegraphics[width=2.9in]{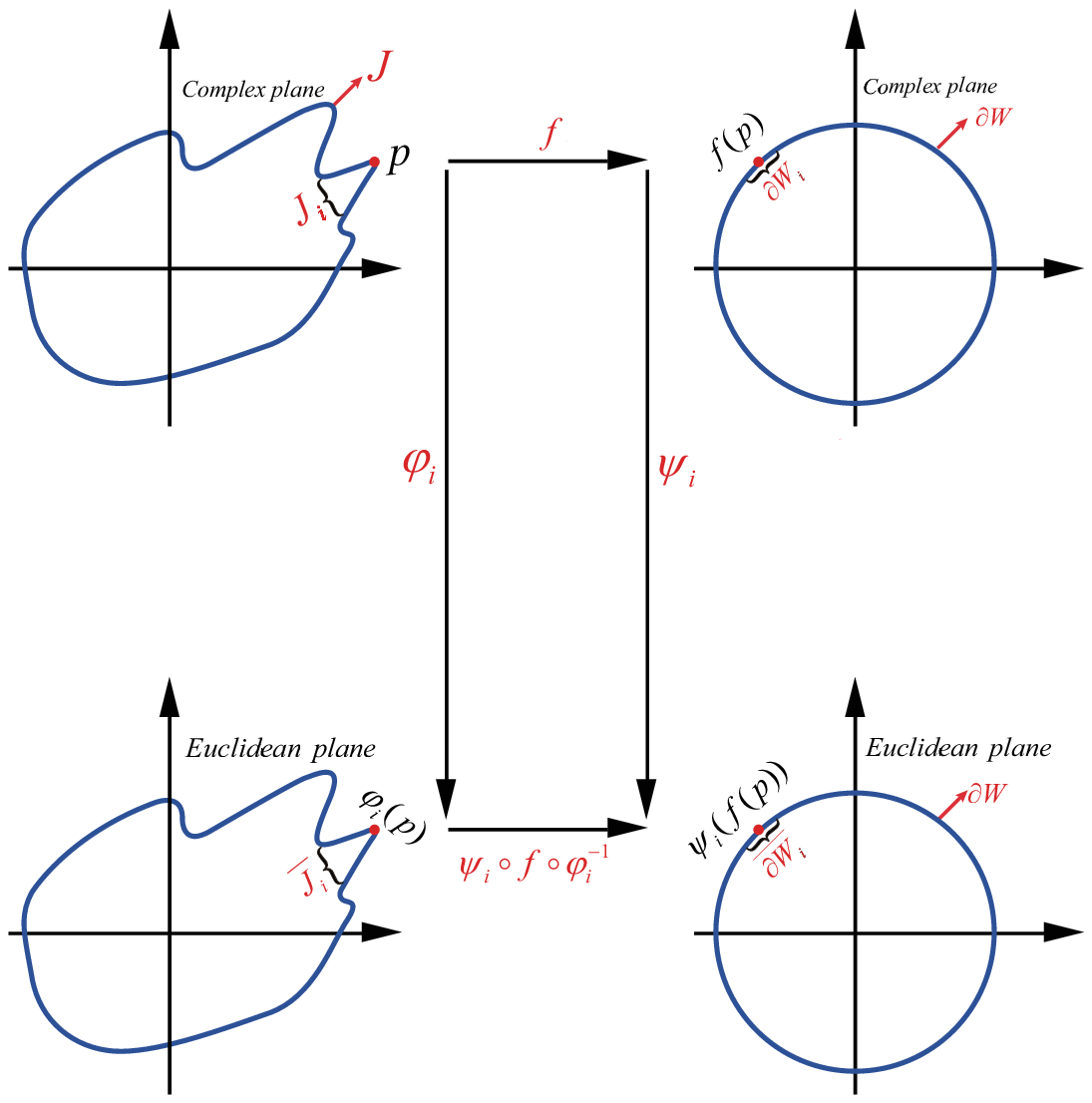}
\caption{Differentiable mapping}
\label{figDM}
\end{figure}

    Before proving the lemma, we expressed the simple closed curve as its complex form. Here, we will discuss manifolds from the complex plane to the Euclidean plane.
   Arbitrary simple closed curve $J$ and unit circle $\partial W$ in complex plane are one-dimensional differentiable manifolds.
    For the mapping $f:J \to \partial W$, $J $ and $\partial W $ can, respectively, obtain a coordinate in the Euclidean plane by the chart $\left(J_{i}, \varphi_{i}\right)$(belong to atlas $\left\{\left(J_{i}, \varphi_{i}\right)\right\} )$ and$\left(\partial W_{i}, \psi_{i}\right)$(belong to atlas $\left\{\left(\partial W_{i}, \psi_{i}\right)\right\}$).
    For $\forall p \in {J}$, the arbitrary selected coordinate $({J_i},{\varphi _i})$ satisfies $p \in {J_i}$ and ${\varphi _i}(p) \in \overline {{J_i}} $, the coordinate $({\partial W_i},{\psi _i})$ satisfies $f(p) \in \partial W_{i}$ and ${\psi _i}(f(p)) \in \overline{\partial W_{i}} $.
    Then, the mapping of domains of Euclidean spaces ${\psi _i} \circ f \circ \varphi _i^{ - 1}$ , which is defined in a neighborhood of the point $\varphi_{i}(p)$, must be differentiable.

    The geometric interpretation of complex number is noted: The one to one correspondence is established between each complex number and ordered pair$(x,y)$ which is a point in plane cartesian coordinate system.
    So the $2$-dimensional Euclidean space can be identified with the $1$-dimensional complex linear space\cite{FomenkoMishchenko-2009-p64}. The specific process is as follows:
  Set $M$ be a $1$-dimensional manifold in complex plane and its atlas is $\left\{U_{\beta}\right\}$. And set $\varphi_{\beta} : U_{\beta} \rightarrow V_{\beta} \subset R^{2}$ be a coordinate homeomorphic mapping. Then, a point $z$ yield two real coordinates $(x, y)$, $z= x + iy $.

    And then, the two coordinate functions ${x_\beta }(p)$ and ${y_\beta }(p)$ in the chart $U_{\beta}$ transform into one complex-valued function ${z_\beta } = {x_\beta }(p) + i{y_\beta }(p)$, ${x_\beta }(p)$ and ${y_\beta }(p)$  are called the complex coordinates of a point in the chart $U_{\beta}$.
   So the coordinate homeomorphic mapping $\psi_{i}$ and $\varphi_{i}$ can be expressed as identity mapping from complex plane to two-dimensional Euclidean space.
   Then, the differentiable mapping ${\psi _i} \circ f \circ \varphi _i^{ - 1}$ can be rewritten as $f$ which is differentiable at point $p$.
   For the arbitrariness of point $p$, we derive that the $f$ is a differentiable mapping from any simple closed curve $J$ to the unit circle $\partial W$.
    Similarly, the inverse mapping $f^{-1}$ is also differentiable.

    Finally, that one arbitrary simple closed curve in plane is diffeomorphic to unit circle is obtained.

    In addition, the lemma can also be proved by the $1$-dimension compact manifold classification theorem\cite{AudinDamian-2014-p49}. $\square$

  \end{description}
\end{description}

\section*{Appendix B}
\appendix
\renewcommand{\theequation}{\thesection B.\arabic{equation}}
In this appendix, we will introduce the other two types of the polar coordinate representation of a planar dynamical system with limit cycle, and discuss wether it is appropriate to use (\ref{206}) to represent the polar coordinate form of a system with limit cycle.

Here, a necessary definition is given.

\textbf{Definition B.1(Equivalent system\cite{MaZhou-2013-p25})}
\emph{If the trajectories (including singular points) of the two autonomous systems are exactly the same (the directions can be different), these two systems are called equivalent.
}

What's more, the reference \cite{MaZhou-2013-p25} has the following statement: For autonomous systems, we mainly study the behavior of their trajectories. Therefore, studying a given autonomous system is the same as studying its equivalent system.

There are the other two types of the polar coordinate representation of a planar dynamical system with limit cycle
\begin{equation}\label{1j6000}
\left\{ {\begin{array}{*{20}{l}}
{\dot r = {\Upsilon _0}(r){\Upsilon _1}(\theta )}\\
{\dot \theta  = \psi (r,\theta )}
\end{array}} \right.,
\end{equation}
and
\begin{equation}\label{1j6001}
\left\{ {\begin{array}{*{20}{l}}
{\dot r = {\Upsilon _0}(r){\Upsilon _2}(r,\theta )}\\
{\dot \theta  = \psi (r,\theta )}
\end{array}} \right.,
\end{equation}
here, we set $r=r^*>0$ to be one limit cycle of system (\ref{1j6000}) and (\ref{1j6001}).

When the planar dynamical system can be transformed to (\ref{1j6000}), we have ${\Upsilon _0}(r^*) \equiv 0$, ${\Upsilon _1}(\theta ) \not \equiv 0$. So the function ${\Upsilon _1}(\theta )$ does not affect the position or size of the limit cycle.

\begin{figure}[h]
\centering
\subfloat[The phase diagram of system (\ref{1a41a}).]{%
\resizebox*{5.9cm}{!}{\includegraphics{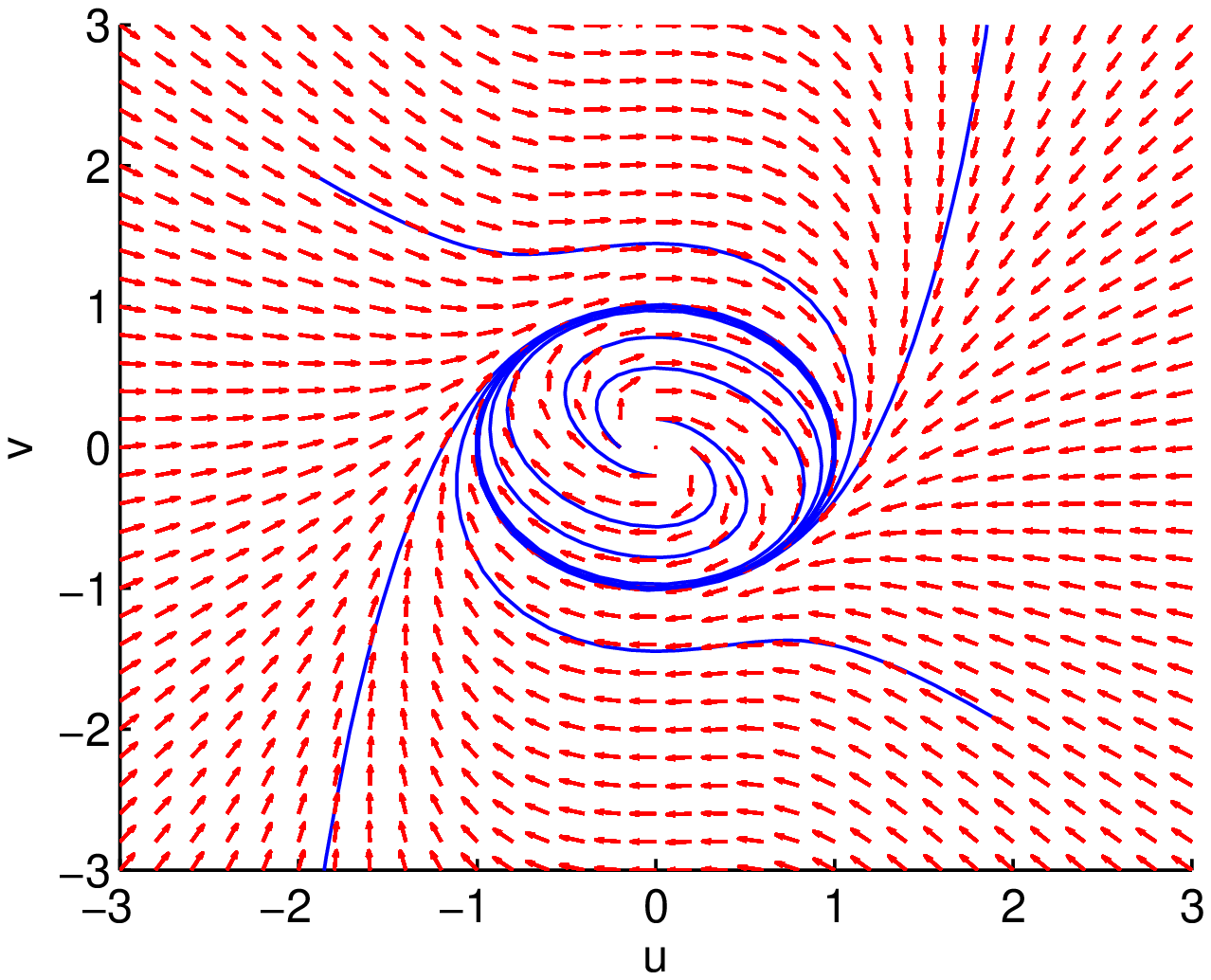}}}\hspace{5pt}
\subfloat[The phase diagram of system (\ref{181}).]{%
\resizebox*{5.9cm}{!}{\includegraphics{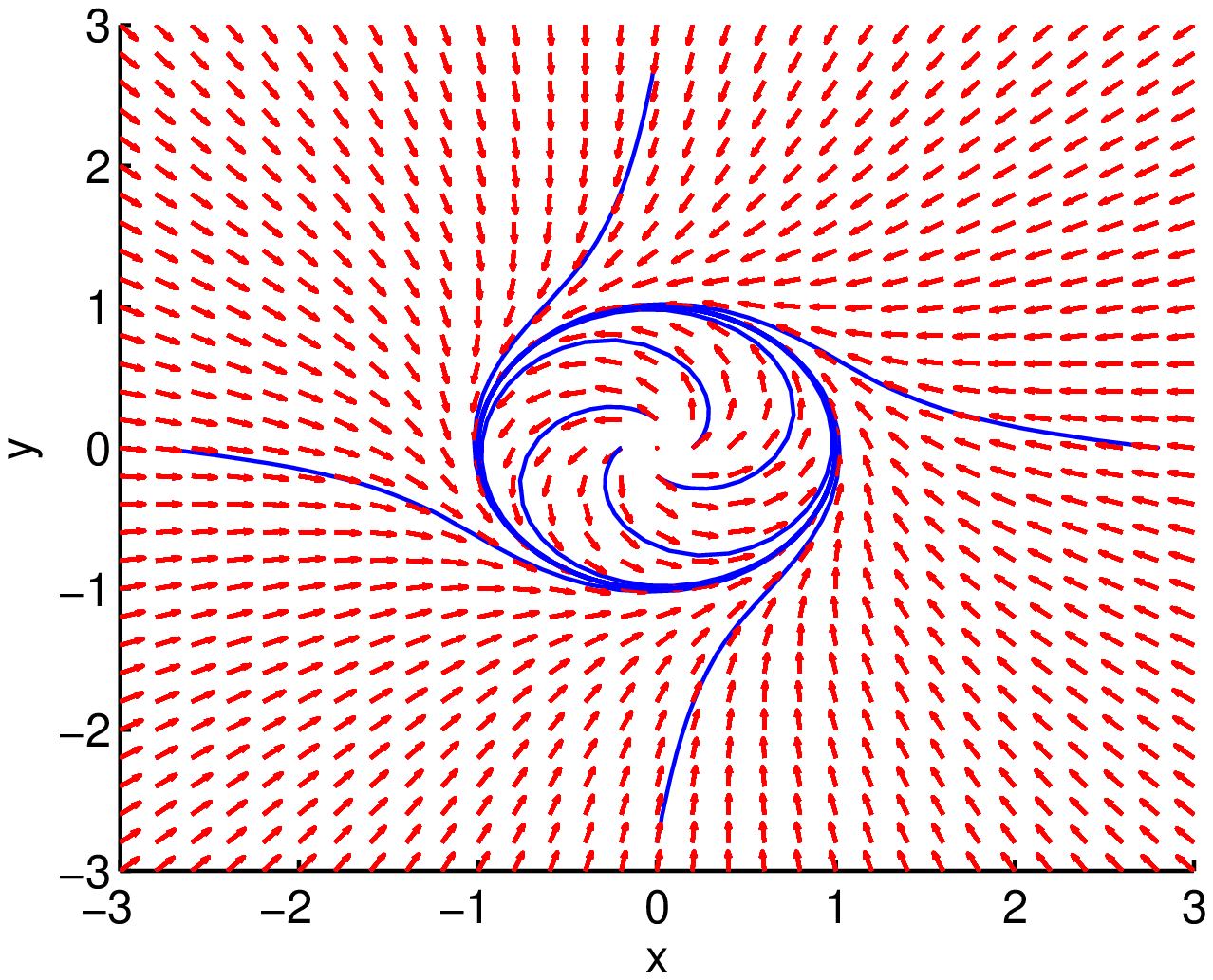}}}
\caption{The phase diagram of system (\ref{1a41a}) and (\ref{181}).}
\label{AppendixBFig1}
\end{figure}

The Figure. \ref{AppendixBFig1} shows the phase diagrams of the special case (\ref{1a41a}) of (\ref{1j6000}) and the special case (\ref{181}) of (\ref{206}).

Then, by the phase diagrams in Figure. \ref{AppendixBFig1} and the Definition B.1, we can obtain the system (\ref{1a41a}) and (\ref{181})  are equivalent systems to each other.
Furthermore, it obtains that the system (\ref{1j6000}) and (\ref{206}) are equivalent systems to each other, too.

When the planar dynamical system can be transformed to (\ref{1j6001}), a simple example \cite{Giesl-2004-p643} is given
\begin{equation}\label{1j1065}
\left\{ {\begin{array}{*{20}{l}}
{\dot x = x\left( {1 - {x^2} - {y^2}} \right)\left( {x + \frac{1}{2}} \right) - y}\\
{\dot y = y\left( {1 - {x^2} - {y^2}} \right)\left( {x + \frac{1}{2}} \right) + x}
\end{array}} \right..
\end{equation}
By the polar transformation, the (\ref{1j1065}) can be transformed into
\begin{equation}\label{1j1066}
\left\{ \begin{array}{l}
\frac{{dr}}{{dt}} = r(1 - {r^2})(r\cos \theta  + \frac{1}{2})\\
\frac{{d\theta }}{{dt}} = 1
\end{array} \right.,
\end{equation}
here, ${\Upsilon _0}(r) = r(1 - {r^2}),{\Upsilon _2}(r,\theta ) = r\cos \theta  + \frac{1}{2}$.

It is easy to find ${\Upsilon _0}(1) \equiv 0,{\Upsilon _2}(r,\theta ) \not \equiv 0$. So the system (\ref{1j1066}) takes the unit circle as the limit cycle, but its radial coordinate system have cross terms ${\Upsilon _2}(r,\theta )$.

\begin{figure}[h]
\centering
\subfloat[The phase diagram of system (\ref{1j1066}).]{%
\resizebox*{5.9cm}{!}{\includegraphics{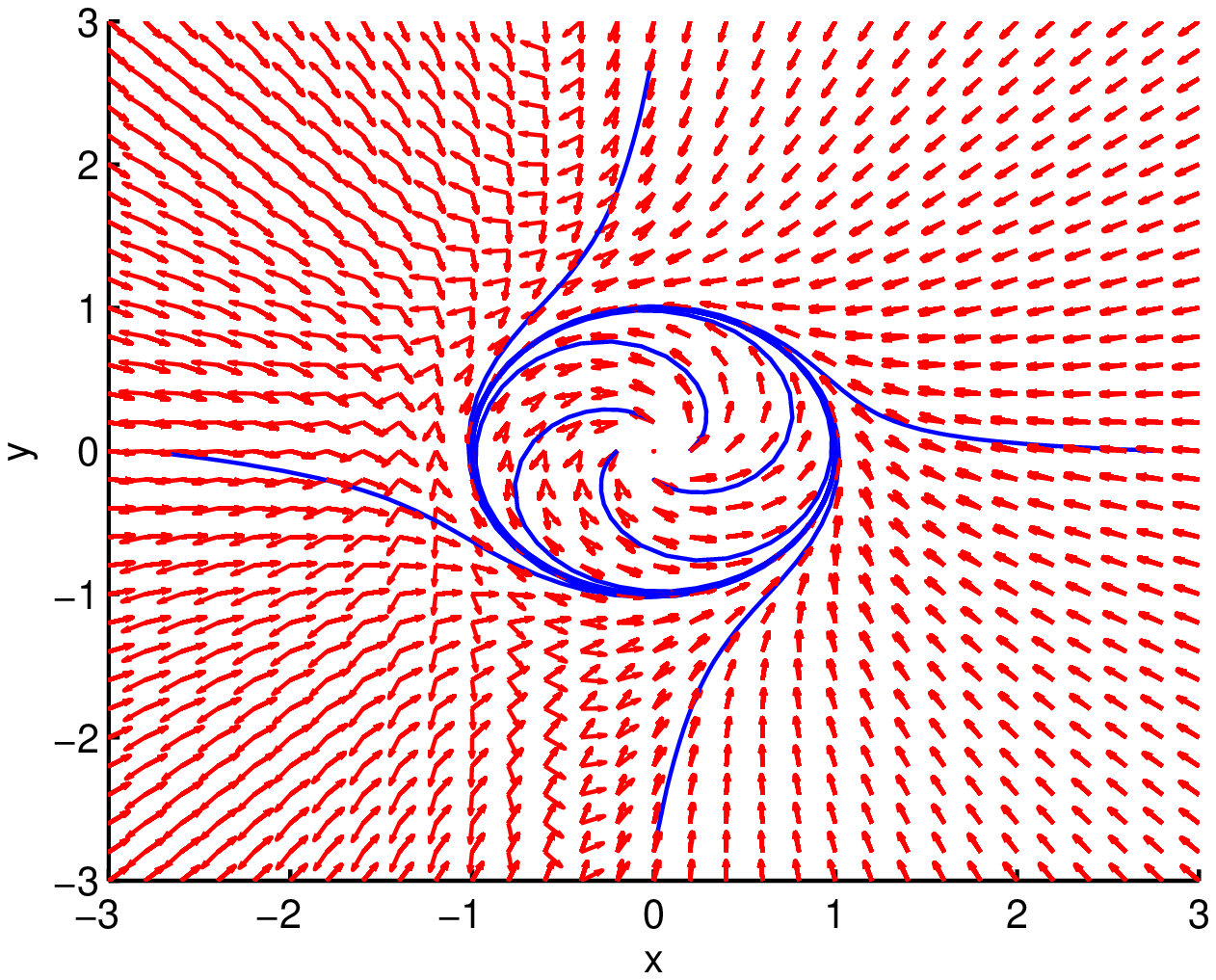}}}\hspace{5pt}
\subfloat[The phase diagram of system (\ref{181}).]{%
\resizebox*{5.9cm}{!}{\includegraphics{ex2a.eps}}}
\caption{The phase diagram of system (\ref{1j1066}) and (\ref{181}).}
\label{AppendixBFig}
\end{figure}

What's more, the Figure. \ref{AppendixBFig} shows the phase diagrams of the special case (\ref{1j1066}) of (\ref{1j6001}) and the special case (\ref{181}) of (\ref{206}).

Then, by the phase diagrams in Figure. \ref{AppendixBFig} and the Definition B.1, we can obtain the system (\ref{1j1066}) and (\ref{181}) are equivalent systems to each other.
Furthermore, it obtains that the system (\ref{1j6001}) and (\ref{206}) are equivalent systems to each other, too.

To sum up, it is appropriate to use system (\ref{206}) to represent the polar coordinate form of a system with one limit cycle.

\begin{backmatter}
\section*{Acknowledgements}
Thanks members in the Institute of Systems Science of Shanghai University for discussion, particularly, with Xin-Jian Xu.
We are also grateful to Wen-Qing Hu of Missouri University of Science and Technology for many useful suggestions.

\section*{Funding}
This work was supported in part by the Natural Science Foundation of China No. NSFC91329301 and No. NSFC9152930016; and by the grants from the State Key Laboratory of Oncogenes and Related Genes (No. 90-10-11).

\section*{Availability of data and materials}
The data sets used or analysed during the current study are available from the corresponding author on reasonable
request.

\section*{Competing interests}
  The authors declare that they have no competing interests.

\section*{Author's contributions}
All authors read and approved the final manuscript.


\bibliographystyle{bmc-mathphys} 
\bibliography{bmc_article}      

\end{backmatter}
\end{document}